\documentclass[twoside,11pt,a4paper,leqno]{article}
\usepackage{amssymb,amsmath,amscd,euscript,verbatim,array}
\usepackage[center,pagestyles]{titlesec}
\usepackage{geometry}
\usepackage{anysize}
\usepackage{fancyhdr}
\usepackage{indentfirst}
\usepackage{graphicx}
\usepackage{color}
\usepackage{ifpdf}
\usepackage{mathrsfs}
\usepackage{cite}
\usepackage[numbers,sort&compress]{natbib}
\marginsize{3.5cm}{3.5cm}{2cm}{2cm}

\titleformat{\section}{\centering\normalsize}{\thesection.}{0.5em}{}
\titleformat{\subsection}{\normalsize\bfseries}{\thesubsection.}{0.5em}{}
\titleformat{\subsubsection}{\normalsize\bfseries}{\thesubsubsection.}{0.5em}{}
\newcommand{\N}{\mathbb{N}}

\newcommand{\R}{\mathbb{R}}

\newtheorem{Theorem}{Theorem}[section]
\newtheorem{Definition}[Theorem]{Definition}
\newtheorem{Lemma}[Theorem]{Lemma}
\newtheorem{Exercise}[Theorem]{Exercise}
\newtheorem{Proposition}[Theorem]{Proposition}
\newtheorem{Remark}[Theorem]{Remark}

\newcommand{\gm}{\gamma}

\newcommand{\bthm}{\begin{Theorem}}
\newcommand{\ethm}{\end{Theorem}}
\newcommand{\bpr}{\begin{Proposition}}
\newcommand{\epr}{\end{Proposition}}
\newcommand{\blm}{\begin{Lemma}}
\newcommand{\elm}{\end{Lemma}}
\newcommand{\bex}{\begin{Exercise}}
\newcommand{\eex}{\end{Exercise}}
\newcommand{\be}{\begin{equation}}
\newcommand{\ee}{\end{equation}}
\newcommand{\beal}{\begin{aligned}}
\newcommand{\enal}{\end{aligned}}
\newcommand{\brm}{\begin{Remark}}
\newcommand{\erm}{\end{Remark}}
\newcounter{item}[section]

\newcommand{\Proof}{\textbf{Proof}\hspace{0.3cm}}
\newcommand{\End}{\ensuremath{\hfill{\Box}}\\}
\renewcommand{\title}[1]{\begin{center}\textbf{\large #1}\end{center}}
\renewcommand{\author}[1]{\begin{center}\normalsize #1\end{center}}
\renewcommand{\date}[1]{\begin{center}#1\end{center}}

\makeatletter \@addtoreset{equation}{section}
\makeatother
\begin{document}
\vspace{10pt}
\title{WEAK KAM THEORY FOR GENERAL HAMILTON-JACOBI EQUATIONS III: THE VARIATIONAL PRINCIPLE UNDER OSGOOD CONDITIONS}

\vspace{6pt}
\author{\sc Lin Wang and Jun Yan}

\vspace{10pt} \thispagestyle{plain}

\begin{quote}
\small {\sc Abstract.} We consider the following evolutionary Hamilton-Jacobi equation with initial condition:
\begin{equation*}
\begin{cases}
\partial_tu(x,t)+H(x,u(x,t),\partial_xu(x,t))=0,\\
u(x,0)=\phi(x),
\end{cases}
\end{equation*}
where $\phi(x)\in C(M,\R)$.
Under some assumptions on the convexity  of  $H(x,u,p)$ with respect to $p$ and the Osgood growth of  $H(x,u,p)$ with respect to $u$, we establish an implicitly variational principle and provide  an intrinsic relation between  viscosity solutions and certain minimal characteristics. Moreover, we obtain a representation formula of the viscosity solution of the evolutionary  Hamilton-Jacobi equation.
\end{quote}
\begin{quote}
\end{quote}
\begin{quote}
\end{quote} \vspace{25pt}

\tableofcontents
\newpage
\section{\sc Introduction and main results}
Let $M$ be a closed manifold  and $H$ be a $C^r$ ($r\geq 2$) function called a Hamiltonian. We consider the following Hamilton-Jacobi equation:
\begin{equation}\label{hje}
\partial_tu(x,t)+H(x,u(x,t),\partial_xu(x,t))=0,
\end{equation}with the initial condition
\[u(x,0)=\phi(x),\]where $(x,t)\in M\times[0,T]$, $T$ is a positive constant.
The characteristics of (\ref{hje}) satisfies the following equation:
\begin{equation}\label{hjech}
\begin{cases}
\dot{x}=\frac{\partial H}{\partial p},\\
\dot{p}=-\frac{\partial H}{\partial x}-\frac{\partial H}{\partial u}p,\\
\dot{u}=\frac{\partial H}{\partial p}p-H.
\end{cases}
\end{equation}
To avoid the ambiguity, we denote the solution of (\ref{hjech}) (the characteristics of (\ref{hje})) by $(X(t),U(t),P(t))$.

In 1983, M. Crandall and P. L. Lions introduced a notion of weak solution named viscosity solution for overcoming the lack of uniqueness of the solution due to the crossing of characteristics (see \cite{Ar, CL2}). Owing to the notion itself, the uniqueness of the viscosity solution can be followed from comparison principle (see \cite{Ba1,Ba2,CEL,CHL1,CHL2,CL2} for instance). However, the nondecreasing property of $H(x,u,p)$ with respect to $u$ was necessary to achieve the uniqueness of the viscosity solution. More generally, it was required that for certain $\gm\in\R$, $H(x,u,p)-\gm u$ is nondecreasing with respect to $u$.
During the same period, S. Aubry and J. Mather developed a seminar work so called Aubry-Mather theory on global action minimizing orbits for area-preserving twist maps (see \cite{Au,AD,M0,Mat,M2,M3} for instance). Moreover, it was generalized to positive definite Lagrangian systems with multi-degrees of freedom in \cite{M1}.

There is a close connection between viscosity solutions  and Aubry-Mather theory. Roughly speaking, the global minimizing orbits used in Aubry-Mather theory can be embedded into the characteristic fields of PDEs. The similar ideas were reflected  in pioneering papers \cite{E} and \cite{F2} respectively. In \cite{E}, W. E  was concerned with certain weak solutions of Burgers equation. In \cite{F2}, A. Fathi provided a weak solution named weak KAM solution  and implied that the weak KAM solution is a viscosity solution, which initiated so called weak KAM theory.   Later, it was obtained the equivalence between weak KAM solutions and viscosity solutions for the Hamiltonian $H(x,p)$ without the unknown function $u$ under strict convexity and superlinear growth  with respect to $p$.  Moreover, based on the relations between weak KAM solutions and viscosity solutions, the regularity of global subsolutions was improved (see \cite{Be,FS}). A systematic introduction to weak KAM theory can be found in \cite{F3}.

Due to the lack of the variational principle for more general Hamilton-Jacobi equations, the weak KAM theory had been limited to Hamilton-Jacobi equations without the unknown function $u$ explicitly. In \cite{SWY}, the authors made an attempt on the Hamilton-Jacobi equation formed as (\ref{hje}) by a dynamical approach and extended Fathi's weak KAM theory to slightly general Hamilton-Jacobi equations under the  monotonicity (non-decreasing) which is also referred as ``proper" condition. In particular, the convergence of the viscosity solutions of evolutionary equations were obtained. In \cite{WY1},  A more general weak KAM theory was established without the monotonicity assumption. Unfortunately, the convergence of the viscosity solutions do not holds generally. In both  \cite{SWY} and \cite{WY1}, the assumption on uniformly Lipschitz  of $H$ with respect to $u$ was needed, which still built a barrier of the scope of the weak KAM theory.

In this paper, both of the monotonicity and uniformly Lipschitz above are replaced by a more general Osgood growth assumption, which is called ``Osgood growth" (see (H4)), which makes a further step to  enlarge the scope of the weak KAM theory. More precisely, we establish a variational principle and provide  an intrinsic relation between  viscosity solutions and certain minimal characteristics. Moreover, we obtain a representation formula of the viscosity solution of (\ref{hje}). We are concerned with the viscosity solutions with finite time in this context. The large time behavior of the solutions can be discussed based on similar arguments as \cite{WY1}. Precisely speaking, we are concerned with a  $C^r$ ($r\geq 2$) Hamiltonian $H(x,u,p)$ satisfying the following conditions:
\begin{itemize}
\item [\textbf{(H1)}] \textbf{Positive Definiteness}: $H(x,u,p)$ is strictly convex with respect to $p$;
\item [\textbf{(H2)}] \textbf{Superlinearity in the Fibers}: For every compact set $I$ and any $u\in I$, $H(x,u,p)$ is uniformly superlinear growth with respect  to $p$;
\item [\textbf{(H3)}] \textbf{Completeness of the Flow}: The flows of (\ref{hjech}) generated by $H(x,u,p)$ are complete.
\item [\textbf{(H4)}] \textbf{Osgood Growth}: For every compact set $K$ and any $(x,p)\in K\subset T^*M$, there exists a continuous function $f_K(u)$ defined on $[0,+\infty)$
    with the divergent integral $\int_0^\infty\frac{1}{f_K(u)}du$  such that
    \[H(x,|u|,p)\geq \langle \frac{\partial H}{\partial p},p\rangle-f_K(|u|).\]
\end{itemize}
It is easy to see that $\int_0^\infty\frac{1}{f_K(u)}du$ is divergent if and only if the flow generated by $\dot{u}=f_K(u)$ is complete. (H4)  can be referred as Osgood condition (see \cite{os}).

We use $\mathcal{L}: T^*M\rightarrow TM$ to denote  the Legendre transformation. Let
$\bar{\mathcal{L}}:=(\mathcal{L}, Id)$, where $Id$ denotes the identity map from $\R$ to $\R$. Then $\bar{\mathcal{L}}$ denote a diffeomorphism from $T^*M\times\R$ to $TM\times\R$. By $\bar{\mathcal{L}}$,
the Lagrangian $L(x,u, \dot{x})$ associated to $H(x,u,p)$ can be denoted by
\[L(x,u, \dot{x}):=\sup_p\{\langle \dot{x},p\rangle-H(x,u,p)\}.\]
Let $\Psi_t$ denote the flows of (\ref{hjech}) generated by $H(x,u,p)$. The flows generated by $L(x,u,\dot{x})$
can be denoted by $\Phi_t:=\bar{\mathcal{L}}\circ\Psi_t\circ\bar{\mathcal{L}}^{-1}$. Based on
(H1)-(H4), it follows from $\bar{\mathcal{L}}$ that the Lagrangian $L(x,u,\dot{x})$ satisfies:
\begin{itemize}
\item [\textbf{(L1)}] \textbf{Positive Definiteness}: $L(x,u,\dot{x})$ is strictly convex with respect  to $\dot{x}$;
\item [\textbf{(L2)}] \textbf{Superlinearity in the Fibers}: For every compact set $I$ and any $u\in I$, $L(x,u,\dot{x})$  is uniformly superlinear growth with respect  to $\dot{x}$;
\item [\textbf{(L3)}] \textbf{Completeness of the Flow}: The flows generated by $L(x,u,\dot{x})$  are complete.
\item [\textbf{(L4)}] \textbf{Osgood Growth}: For every compact set $K$ and any $(x,\dot{x})\in K\subset TM$, there exists a continuous function $f_K(u)$ defined on $[0,+\infty)$
    with the divergent integral $\int_0^\infty\frac{1}{f_K(u)}du$  such that
    \[L(x,|u|,\dot{x})\leq f_K(|u|).\]
\end{itemize}
It is easy to see that (L4) is more general than the  monotonicity (non-increasing)  and Lipschitz  of $L$ with respect to $u$, for which $f_K(u)$ is corresponding to a constant function and an affine function respectively.

If a Hamiltonian $H(x,u,p)$ satisfies (H1)-(H4), then we obtain the following theorem:
\begin{Theorem}\label{two}
For given $x_0, x\in M$, $u_0\in\R$ and $t\in (0, T]$, there exists a unique $h_{x_0,u_0}(x,t)$ satisfying
\begin{equation}
h_{x_0,u_0}(x,t)=u_0+\inf_{\substack{\gm(t)=x \\  \gm(0)=x_0} }\int_0^tL(\gm(\tau),h_{x_0,u_0}(\gm(\tau),\tau),\dot{\gm}(\tau))d\tau,
\end{equation}where the infimums are taken among the continuous and piecewise $C^1$ curves. In particular, the infimums are attained at the characteristics of (\ref{hje}).
 Moreover, let $\mathcal{S}_{x_0,u_0}^x$ denote the set of characteristics $(X(t),U(t),P(t))$ satisfying $X(0)=x_0$, $X(t)=x$ and $U(0)=u_0$, then we have
\begin{equation}
h_{x_0,u_0}(x,t)=\inf\left\{U(t):(X(t),U(t),P(t))\in \mathcal{S}_{x_0,u_0}^x\right\}.
\end{equation}
\end{Theorem}

Theorem \ref{two} provides a general variational principle, which builds a bridge between Hamilton-Jacobi equations under (H1)-(H4) and Hamiltonian dynamical systems. As an application, we will obtain a dynamical representation of the viscosity solution of (\ref{hje}) in the following. More precisely, we have the following theorem:
\begin{Theorem}\label{one}
There exists a unique  viscosity solution $u(x,t)$ of (\ref{hje}) with initial condition $u(x,0)=\phi(x)$. Moreover, $u(x,t)$  can be represented as
\begin{equation}
u(x,t)=\inf_{y\in M}h_{y,\phi(y)}(x,t).
\end{equation}
\end{Theorem}

Theorem \ref{two} and Theorem \ref{one} implies the following theorem directly:
\begin{Theorem}\label{three}
For $(x,t)\in M\times [0,T]$, the  viscosity solution $u(x,t)$ of (\ref{hje}) with initial condition $u(x,0)=\phi(x)$ is determined by the minimal characteristic curve. More precisely, we have
\begin{equation}
u(x,t)=\inf_{y\in M}\inf\left\{U(t):(X(t),U(t),P(t))\in \mathcal{S}_{y,\phi(y)}^x\right\},
\end{equation}
where $\mathcal{S}_{y,\phi(y)}^x$ denotes the set of characteristics $(X(t),U(t),P(t))$ satisfying $X(0)=y$, $X(t)=x$ and $U(0)=\phi(y)$.
\end{Theorem}

A similar result corresponding to the viscosity solutions of Hamilton-Jacobi equations without the unknown function $u$ was well known (see Theorem 6.4.6 in \cite{CS} for instance). Theorem \ref{three} implies the relation between the viscosity solutions and  the minimal characteristics still holds for more general Hamilton-Jacobi equations. Roughly speaking, the notion of viscosity solution was invented to avoid the lack of uniqueness owing to the crossing of characteristics.  Based on Theorem \ref{three}, the reason why  the notion of viscosity solution results in the fact without crossing is that the properties of viscosity solutions are determined by certain minimal characteristics.


\section{\sc Preliminaries}
In this section, we recall the definitions  and some properties of the viscosity solution of (\ref{hje})  (see \cite{CEL,CL2,F3}). In addition, we provide some aspects of Mather-Fathi theory for  the sake of completeness.

\subsection{Viscosity solutions and semiconcavity}
We introduce the notions of upper and lower differentials (see \cite{CS,F3} for instance).
\begin{Definition}
Let $u:U\rightarrow\R$ be a function defined on the open subset $U$ of $\R^k$, then the linear form $\theta$ is called a upper  differential of $u$ at $x_0$ if
\[\limsup_{x\rightarrow x_0}\frac{u(x)-u(x_0)-\theta(x-x_0)}{|x-x_0|}\leq 0;\]
In the same way, $\theta$ is called a lower  differential of $u$ at $x_0$ if
\[\liminf_{x\rightarrow x_0}\frac{u(x)-u(x_0)-\theta(x-x_0)}{|x-x_0|}\geq 0.\]
\end{Definition}

Let $U$ be an open convex subset of $\R^k$ and let $u:U\rightarrow\R$ be a function. $u$ is called a semiconcave function if there exists a finite constant $K$ and for each $x\in U$ there exists a linear form $\theta_x:\R^k\rightarrow\R$ such that for any $y\in U$
\begin{equation}
u(y)-u(x)\leq \theta_x(y-x)+K|y-x|^2.
\end{equation}
For the sake of simplicity, we only consider the semiconcave functions with linear modulus defined as above. See \cite{CS} for a more general definition. In this context, the notion ``semiconcave" means ``semiconcave with a linear modulus".
\begin{Definition}
A function $u:M\rightarrow\R$ defined on the $C^r$ ($r\geq 2$) differential $k$-dimensional manifold $M$ is locally semiconcave if for each $x\in M$ there exists a $C^r$ ($r\geq 2$) coordinate chart $\psi:U\rightarrow\R^k$ with $x\in U$ such that $u\circ \psi^{-1}:U\rightarrow\R$ is semiconcave.
\end{Definition}

Following from \cite{CEL,CL2,F3},
 a viscosity solution of (\ref{hje}) can be defined as follows:
\begin{Definition}\label{visco}
Let $V$ be an open subset  $V\subset M$,
\begin{itemize}
\item [(i)] A function $u:V\times[0,T]\rightarrow \R$ is a subsolution of  (\ref{hje}), if for every $C^1$ function $\phi:V\times[0,T]\rightarrow\R$ and every point $(x_0,t_0)\in V\times[0,T]$ such that $u-\phi$ has a maximum at $(x_0,t_0)$, we have
\begin{equation}
\partial_t\phi(x_0,t_0)+H(x_0,u(x_0,t_0),\partial_x\phi(x_0,t_0))\leq 0;
\end{equation}
\item [(ii)] A function $u:V\times[0,T]\rightarrow \R$ is a supersolution of  (\ref{hje}), if for every $C^1$ function $\psi:V\times[0,T]\rightarrow\R$ and every point $(x_0,t_0)\in V\times[0,T]$ such that $u-\psi$ has a minimum at $(x_0,t_0)$, we have
\begin{equation}
\partial_t\psi(x_0,t_0)+H(x_0,u(x_0,t_0),\partial_x\psi(x_0,t_0))\geq 0;
\end{equation}
\item [(iii)] A function $u:V\times[0,T]\rightarrow \R$ is a viscosity solution of  (\ref{hje}) on the open subset  $V\subset M$, if it is both a subsolution and a supersolution.
\end{itemize}
\end{Definition}

Based on \cite{CS} (Theorem 5.3.1. and Theorem 5.3.6), we have the following results.
\begin{Proposition}\label{semi}
Considering the Hamilton-Jacobi equations
\[H(x,u,\partial_xu)=0 \quad\text{and}\quad \partial_tu+H(x,u,\partial_xu)=0.\]
 where $H\in C^r(T^*M\times\R,\R)$ ($r\geq 2$), we have the following properties.
\begin{itemize}
\item [(a)] Let $u$ be a semiconcave function satisfying the equation almost everywhere. If $H(x,u,p)$ is  convex with respect to $p$, then $u$ is a viscosity solution of the equation;
\item [(b)] Let $u$ be a Lipschitz viscosity solution of the equation. If $H(x,u,p)$ is strictly convex with respect to $p$, then $u$ is locally semiconcave on $M$.
\end{itemize}
\end{Proposition}

\subsection{The minimal action and the  fundamental solution}
Let $L:TM\rightarrow\R$ be a Tonelli Lagrangian. We define the function $h_t:M\times M\rightarrow\R$ by
\begin{equation}\label{mather}
h_t(x,y)=\inf_{\substack{\gm(0)=x\\ \gm(t)=y}}\int_0^t L(\gm(\tau),\dot{\gm}(\tau))d\tau,
\end{equation}
where the infimums are taken among the absolutely continuous curves $\gm:[0,t]\rightarrow M$.
By Tonelli theorem (see \cite{F3,M1}), the infimums in (\ref{mather}) can be achived. Let $\bar{\gm}$ be an absolutely continuous curve with $\bar{\gm}(0)=x$ and $\bar{\gm}(t)=y$ such that the infinmum is achieved at $\bar{\gm}$. Then $\bar{\gm}$ is called a minimal curve. By \cite{M1}, the minimal curves satisfy the Euler-Lagrange equation generated by $L$. The quantity $h_t(x,y)$ is called a minimal  action.
From the definition of $h_t(x,y)$, it follows that for each $x,y,z\in M$ and each $t,t'>0$, we have
 \begin{equation}\label{yan}
h_{t+t'}(x,z)\leq h_{t}(x,y)+h_{t'}(y,z).
\end{equation}
  In particular,  we have
\begin{equation}\label{wang}
h_{t+t'}(x,y)=h_{t}(x,\bar{\gm}(t))+h_{t'}(\bar{\gm}(t),y),
\end{equation}
where $\bar{\gm}$ is  a minimal curve with $\bar{\gm}(0)=x$ and $\bar{\gm}(t+t')=y$.

Consider the following Hamilton-Jacobi equation:
 \begin{equation}\label{fathitt}
\begin{cases}
\partial_tu(x,t)+H(x,  \partial_xu(x,t))=0,\\
u(x,0)=\phi(x),
\end{cases}
\end{equation}
where $\phi(x)\in C(M)$.
By \cite{F3}, a viscosity solution of (\ref{fathitt}) can be represented as
 \begin{equation}\label{infc}
u(x,t):=\inf_{y\in M}\left\{\phi(y)+h^t(y,x)\right\}.
\end{equation}
The right side of (\ref{infc}) is also called inf-convolution of $\phi$, due to the formal analogy with the usual convolution (see \cite{CS}). Moreover, the minimal action $h^t(y,x)$ can be viewed as a fundamental solution of (\ref{fathitt}) (see \cite{GT}).

The following conception is crucial in our context.
\begin{Definition}\label{cali}
For $u(x,t)\in C(M\times[0,T],\R)$, a curve $\gm:I\rightarrow M$ is called a calibrated curve of $u$ if for every $t_1, t_2\in I$ with $0\leq t_1< t_2$, we have
\[u(\gm(t_2),t_2)=u(\gm(t_1),t_1)+\int_{t_1}^{t_2}L(\gm(\tau),u(\gm(\tau),\tau),\dot{\gm}(\tau))d\tau.\]
\end{Definition}

We are devoted to detecting the viscosity solution of (\ref{hje}) from a dynamical view.  For given $x_0, x\in M$, $u_0\in\R$ and $t\in (0, T]$, we define formally:
\begin{equation}\label{u}
h_{x_0,u_0}(x,t)=u_0+\inf_{\substack{\gm(t)=x \\  \gm(0)=x_0} }\int_0^tL(\gm(\tau),h_{x_0,u_0}(\gm(\tau),\tau),\dot{\gm}(\tau))d\tau,
\end{equation}where the infimums are taken among the continuous and piecewise $C^1$ curves. It is easy to see that the cure achieving the infimum in the right side of (\ref{u}) is a calibrated curve of $h_{x_0,u_0}(x,t)$. To fix the notions, we call $h_{x_0,u_0}(x,t)$ the fundamental solution of (\ref{hje}).
  In next section, we will show the well-posedness of $h_{x_0,u_0}(x,t)$ under the assumptions (L1)-(L4).

\section{\sc Variational principle}
In this section, we are devoted to proving Theorem \ref{two}. First of all, let us recall the variational principle under uniformly Lipschitz conditions based on \cite{WY1}. Precisely speaking, the following assumption is added.
\begin{itemize}
\item [\textbf{(L4')}]  \textbf{Uniform Lipschitz}: $L(x,u,\dot{x})$ is uniformly Lipschitz with respect to $u$.
\end{itemize}
\begin{Theorem}\label{two1}
Let $L$ satisfy (L1), (L2), (L3) and (L4'). For given $x_0, x\in M$, $u_0\in\R$ and $t\in (0, T]$, there exists a unique $h_{x_0,u_0}(x,t)$ satisfying
\begin{equation}
h_{x_0,u_0}(x,t)=u_0+\inf_{\substack{\gm(t)=x \\  \gm(0)=x_0} }\int_0^tL(\gm(\tau),h_{x_0,u_0}(\gm(\tau),\tau),\dot{\gm}(\tau))d\tau,
\end{equation}where the infimums are taken among the absolutely continuous curves $\gm:[0,t]\rightarrow M$. In particular, the infimums are attained at the characteristics of (\ref{hje}).
 Moreover, let $\mathcal{S}_{x_0,u_0}^x$ denote the set of characteristics $(X(t),U(t),P(t))$ satisfying $X(0)=x_0$, $X(t)=x$ and $U(0)=u_0$, then we have
\begin{equation}
h_{x_0,u_0}(x,t)=\inf\left\{U(t):(X(t),U(t),P(t))\in \mathcal{S}_{x_0,u_0}^x\right\}.
\end{equation}
\end{Theorem}

The proof of Theorem \ref{two1} was proved in \cite{WY1}. We omit it here for the consistency of the context. In the following, we will focus on the relaxation of the assumption from (L4') to (L4).
 The proof will be divided into three steps. In the first step, the Lagrangian function $L(x,u,\dot{x})$  will be truncated and it will be proved that $h_{x_0,u_0}(x,t)$ is independent of the truncated part of $L(x,u,\dot{x})$.  In the second step, combing with  the completeness of flow (L3), it will be showed that the additional uniformly Lipschitz assumption on $L$ is also not necessary. Finally,  the proof of Theorem \ref{two} will be completed by an argument on a limit process.

\subsection{Truncation of the Lagrangian function}
In this step, we will provide a priori estimate of $h_{x_0,u_0}(\gm(s),s)$, where $\gm(s):[0,t]\rightarrow M$ is a calibrated curve connecting $x_0$ and $x$.

 For the simplicity of notations, we denote
  \begin{equation}\label{199}
  V(x,u,\dot{x}):=L(x,u,\dot{x})-L(x,0,\dot{x}).
   \end{equation}We choose a $C^2$ function $\rho(u)$ such that $\rho_R(u)=1$ for $|u|\leq R$, $\rho_R(u)=0$ for $|u|>R+1$, otherwise $0<\rho_R(u)<1$. Without loss of generality, one can require $|\rho'_R(u)|<2$. Moreover,
we denote
\begin{equation}\label{200}
V_R(x,u,\dot{x}):=\rho_R(u)V(x,u,\dot{x}).
\end{equation}From the uniformly Lipschitz continuity of $V(x,u,\dot{x})$ with respect to $u$, it follows that $V_R(x,u,\dot{x})$ is also uniformly Lipschitz with respect to $u$. Without ambiguity, we  denote the Lipschitz constant of $V_R(x,u,\dot{x})$ by $\lambda_R$.
From (\ref{199}) and (\ref{200}), it follows that  $V_R(x,0,\dot{x})=0$. Moreover, we have
\begin{equation}\label{br}
|V_R(x,u,\dot{x})|=|V_R(x,u,\dot{x})-V_R(x,0,\dot{x})|\leq \lambda_R|u|\leq \lambda_R (R+1).
\end{equation}
 Let
\begin{equation}\label{lr}
L_R(x,u,\dot{x})=L(x,0,\dot{x})+V_R(x,u,\dot{x}).
\end{equation}
We omit the subscripts $x_0,u_0$ of $h_{x_0,u_0}(x,t)$ for simplicity.   Based on Theorem \ref{two1}, we have that there exists  a function denoted by $h_R(x,t)$ such that
\begin{equation*}
h_R(x,t)=u_0+\inf_{ \gm(t)=x \atop  \gm(0)=x_0 }\int_0^tL(\gm(\tau),0,\dot{\gm}(\tau))+V_R(\gm(\tau),h_R(\gm(\tau),\tau),\dot{\gm}(\tau))d\tau.
\end{equation*}In addition, the curve achieving the infimum is a calibrated curve. Moreover, we have the following lemma.
\begin{Lemma}\label{A}
For a given $(x,t)\in M\times (0,T]$, let $\gm_R(s):[0,t]\rightarrow M$ be a calibrated curve of $h_R$ satisfying $\gm_R(0)=x_0$ and $\gm_R(t)=x$, then  there exists $A^*$ such that for any $s\in [0,t]$
\begin{equation}\label{astar}
|h_R(\gm_R(s),s)|\leq A^*,
\end{equation}
where $A^*$ is a positive constant only depending on $t$.
\end{Lemma}
\Proof
On one hand, we prove $h_R(\gm_R(s),s)$ is lower bounded. By (L2), it yields that $L(\gm_R(s),0,\dot{\gm}_R(s))$ has a lower bound denoted by $-C_1$, where we use $C_i$ to denote the positive constants independent of $R$. By contradiction, we assume that there exists  subsequences $R_n$ and $s_n$ such that
\begin{equation}\label{81}
h_{R_n}(\gm_{R_n}(s_n),s_n)<-n.
\end{equation}
For the simplicity of notations, we denote $u_R(s):=h_R(\gm_R(s),s)$. Since $u_R(0)=h_R(\gm_R(0),0)=u_0$ for any $R$, then it follows from (\ref{81}) that for $n$ large enough, there exists a time interval $[t_1,t_2]$ such that $u_{R_n}(t_1)=u_0$, $u_{R_n}(t_2)=u_0-1$ and $u_{R_n}(s)\in [u_0-1,u_0]$ for $s\in[t_1,t_2]$. By Theorem \ref{two}, we have $\dot{u}_{R_n}(s)=L(\gm_{R_n}(s),u_{R_n}(s),\dot{\gm}_{R_n}(s))$. Hence,
\begin{equation}\label{822}
-1=\int_{t_1}^{t_2}L(\gm_{R_n}(s),u_{R_n}(s),\dot{\gm}_{R_n}(s))ds.
\end{equation}
From (L2), $\dot{u}_{R_n}(s)$ has a lower bound independent of $R_n$ for $s\in [t_1,t_2]$, which is denoted by $-C_2$. Then we have $t_2-t_1\geq 1/C_2,$ which is independent of $R_n$. From (L2), the formula (\ref{822}) implies there exists a subsequence $\tau_n\in [t_1,t_2]$ such that $|\dot{\gm}_{R_n}(\tau_n)|$ has a upper bound independent of $R_n$. Extracting a subsequence if necessary, it follows from the compactness of $M$ that for $n\rightarrow \infty$
\begin{equation}\label{82}
(\tau_n, \gm_{R_n}(\tau_n),u_{R_n}(\tau_n),\dot{\gm}_{R_n}(\tau_n))\rightarrow (\bar{s}, \bar{x},\bar{u},\bar{v}).
\end{equation}
Let $\Phi_s^{R_n}$ be the flow generated by $L_{R_n}(x,u,\dot{x})$. Then it follows from the completeness of the flow that for any $s\in [t_1,t_2]$, $\Phi_{s}^{R_n}(\bar{x},\bar{u},\bar{v})$ is well defined. Theorem \ref{two} implies $(\gm_{R_n}(s),u_{R_n}(\gm_{R_n}(s),s),\dot{\gm}_{R_n}(s))$ is the flow generated by $L_{R_n}(x,u,\dot{x})$. Based on the construction of $L_{R_n}(x,u,\dot{x})$ (see (\ref{lr})), it yields that for any $s\in [t_1,t_2]$, as $n\rightarrow\infty$
\[\Phi_s^{R_n}(\bar{x},\bar{u},\bar{v})\rightarrow\Phi_s(\bar{x},\bar{u},\bar{v}),\]
where the notation $\rightarrow$ means the convergence in the $C^0$-norm of each component of the flow. Moreover, we have
\begin{align*}
&|(\gm_{R_n}(s_n),u_{R_n}(s_n),\dot{\gm}_{R_n}(s_n))-\Phi_{s_n-\bar{s}}(\bar{x},\bar{u},\bar{v})|\\
\leq &|\Phi_{s_n-\tau_n}^{R_n}(\gm_{R_n}(\tau_n),u_{R_n}(\tau_n),\dot{\gm}_{R_n}(\tau_n))-
\Phi_{s_n-\tau_n}^{R_n}(\bar{x},\bar{u},\bar{v})|\\
&+|\Phi_{s_n-\tau_n}^{R_n}(\bar{x},\bar{u},\bar{v})-\Phi_{s_n-\tau_n}(\bar{x},\bar{u},\bar{v})|\\
&+|\Phi_{s_n-\tau_n}(\bar{x},\bar{u},\bar{v})-\Phi_{s_n-\bar{s}}(\bar{x},\bar{u},\bar{v})|.
\end{align*}
 We consider $(\bar{x},\bar{u},\bar{v})$ as the initial condition of $\Phi_s^{R_n}$ at $s=\bar{s}$. Based on the continuous dependence of solutions of ODEs on initial conditions, it follows that for $n$ large enough,
 \begin{equation}\label{90}
|(\gm_{R_n}(s_n),u_{R_n}(s_n),\dot{\gm}_{R_n}(s_n))-\Phi_{s_n-\bar{s}}(\bar{x},\bar{u},\bar{v})|\leq \epsilon,
\end{equation}
where $\epsilon$ is a small constant independent of $n$. By (L3), it yields that $u_{R_n}(s_n)$ has a bound independent of $n$, which is in contradiction with (\ref{81}).
Therefore,  there exists a constant $A_1$ independent of $R$ such that for $s\in [0,t]$,
\begin{equation}
h_R(\gm_R(s),s)\geq A_1.
\end{equation}

On the other hand, we prove $h_R(\gm_R(s),s)$ is upper bounded. First of all, we estimate $h_R(x,t)$ for a given $(x,t)\in M\times (0,T]$. Without loss of generality, we assume $h_R(x,t)>0$. Let $\bar{\gm}_R(s):[0,t]\rightarrow M$ be a straight line satisfying $\bar{\gm}_R(0)=x_0$ and $\bar{\gm}_R(t)=x$.
 Since $h_R(x,t)>0$,  we have the following dichotomy:
  \begin{itemize}
 \item [(I)]  there exists $s_0\in (0,t)$ such that $h_R(\bar{\gm}_R(s_0),s_0)=0$ and $h_R(\bar{\gm}_R(s),s)\geq 0$ for any $s\in [s_0,t]$;
 \item [(II)] for any  $s\in (0,t)$,  $h_R(\bar{\gm}_R(s),s)> 0$.
 \end{itemize}
It suffices to consider Case (I) and Case (II) can be obtained by a similar argument. Without ambiguity, we denote $h_R(s):=h_R(\bar{\gm}_R(s),s)$ for the sake of the simplicity. Since $\bar{\gm}_R$ is a straight line, then it follows from the compactness that for $s\in [0,t]$, $(\bar{\gm}_R(s),\dot{\bar{\gm}}_R(s))$ is contained in a compact set independent of $R$ denoted by $K$.
 For Case (I),  we have
\begin{align*}
h_R(t)\leq h_R(s_0)+\int_{s_0}^{t}L_R(\bar{\gm}_R(\tau),h_R(\tau),\dot{\bar{\gm}}_R(\tau))d\tau
\leq \int_{s_0}^{t}f_K(h_R(\tau))d\tau,
\end{align*}where the second inequality is from the Osgood growth assumption (L4).
Let $g_R(\tau)$ be a function defined on $[0,t-s_0]$ and satisfy
\begin{equation}
\begin{cases}
\dot{g}_R(\tau)= f_K(g_R(\tau)),\\
g_R(0)=h_R(s_0)=0.
\end{cases}
\end{equation}
Hence, we have
\begin{equation}
\int_0^{g_R}\frac{1}{f_K(g_R)}dg_R=\int_0^sd\tau.
\end{equation}
In particular, we have
\begin{equation}
\int_0^{g_R(t-s_0)}\frac{1}{f_K(g_R)}dg_R=t-s_0,
\end{equation}which together with (L4) yields $g_R(t-s_0)$ has a upper bound independent of $R$. By the comparison theorem of ODEs (see \cite{HS} for instance), we have $h_R(t)\leq g_R(t-s_0)$. Hence, $h_R(\gm_R(t),t)$ has a upper bound independent of $R$ denoted by $u_t$.

Secondly,  Let $\gm_R:[0,t]\rightarrow M$ be a calibrated curve satisfying $\gm_R(0)=x_0$ and $\gm_R(t)=x$. We prove $h_R(\gm_R(s),s)$ is upper bounded for any $s\in (0,t)$.  For the simplicity of notations, we denote $u_R(s):=h_R(\gm_R(s),s)$. In particular, $u_R(t)=u_t$. By contradiction, we assume that there exists  subsequences $R_n$ and $s_n$ such that
\begin{equation}\label{84}
u_{R_n}(s_n)>n.
\end{equation}
 Since $u_R(t)=u_t$ for any $R$, then it follows from (\ref{84}) that for $n$ large enough, there exists a time interval $[t_1,t_2]$ such that $u_{R_n}(t_1)=u_t+1$, $u_{R_n}(t_2)=u_t$ and $u_{R_n}(s)\in [u_t,u_t+1]$ for $s\in[t_1,t_2]$. By Theorem \ref{two}, we have \[\dot{u}_{R_n}(s)=L(\gm_{R_n}(s),u_{R_n}(s),\dot{\gm}_{R_n}(s)).\] Hence,
\begin{equation}\label{85}
-1=\int_{t_1}^{t_2}L(\gm_{R_n}(s),u_{R_n}(s),\dot{\gm}_{R_n}(s))ds.
\end{equation}
From (L2), $\dot{u}_{R_n}(s)$ has a lower bound independent of $R_n$ for $s\in [t_1,t_2]$, which is denoted by $-C_2$. Then we have $t_2-t_1\geq 1/C_2,$ which is independent of $R_n$. By a similar argument as (\ref{90}), it yields that $u_{R_n}(s_n)$ has a bound independent of $n$, which is in contradiction with (\ref{84}).

Therefore, there exists a constant $A_2$ independent of $R$ such that
\begin{equation}
h_R(\gm_R(s),s)\leq A_2.
\end{equation}

So far, it suffices to proof Lemma \ref{A} that we take
\begin{equation}
A^*=\max\{A_1,A_2\}.
\end{equation}
This finishes the proof of Lemma \ref{A}.
\End

\subsection{A priori compactness}
In this step, we will prove a priori estimate of $|\dot{\gm}(s)|$, where $\gm(s):[0,t]\rightarrow M$ is a calibrated curve connecting $x_0$ and $x$.

 We construct a $C^2$ function denoted by $\bar{L}_{R}(x,u,\dot{x})$ satisfying
\begin{equation}
\bar{L}_{R}(x,u,\dot{x}):=\alpha_R (\dot{x})L_R(x,u,\dot{x})+\beta_R(\dot{x}) (\dot{x}^2-R^2)^2,
\end{equation}where $L_R$ is defined as (\ref{lr}) and $\alpha_R (\dot{x})$ is a $C^2$ function satisfying
\begin{equation}
\alpha_R (\dot{x})=\left\{\begin{array}{ll}
\hspace{-0.4em}1,&  |\dot{x}|\leq R+1,\\
\hspace{-0.4em}0,&|\dot{x}|>R+2,\\
\end{array}\right.
\end{equation}
$\beta_R(\dot{x}) $ is defined as \begin{equation}
\beta_R(\dot{x})=\left\{\begin{array}{ll}
\hspace{-0.4em}0,&  |\dot{x}|\leq R,\\
\hspace{-0.4em}\mu_R,&|\dot{x}|>R,\\
\end{array}\right.
\end{equation} where $\mu_R$ is a sufficient large constant.
 Hence, for a given $R_0$, there exists $\lambda_0>0$ (depending on $R_0$) such that for any $u, v\in [-R_0,R_0]$, we have
\begin{equation}
|\bar{L}_{R_0}(x,u,\dot{x})-\bar{L}_{R_0}(x,v,\dot{x})|\leq\lambda_0|u-v|.
\end{equation}
By virtue of the arguments above,  for a given $(x,t) (x\neq x_0, t>0)$, there exists a $C^1$ characteristic curve $\gm_{R_0}:[0,t]\rightarrow M$ such that
\begin{equation*}
h_{R_0}(x,t)=u_0+\int_0^t\bar{L}_{R_0}(\gm_{R_0}(\tau),h_{R_0}(\gm_{R_0}(\tau),\tau),\dot{\gm}_{R_0}(\tau))d\tau,
\end{equation*}where  we drop the subscripts $x_0, u_0$ of $h_{R_0}(x,t)$ for simplicity. First of all, we estimate the initial velocity $\dot{\gm}_R(0)$ for any $R>0$.
\begin{Lemma}\label{xr}
For any $R>0$, $|\dot{\gm}_R(0)|$ has a bound independent of $R$.
\end{Lemma}
\Proof
By contradiction, we assume that there exists a subsequence $R_n$ such that $|\dot{\gm}_{R_n}(0)|\geq n$.
Based on Lemma \ref{A},
\begin{equation}
|h_R(\gm_R(s),s)|\leq A^*,
\end{equation}
where $\gm_R(s):[0,t]\rightarrow M$ is a calibrated curve with $\gm_R(0)=x_0$ and $\gm_R(t)=x$ and $A^*$ is a constant independent of $R$. Hence, for any $R>0$, there exists $s_R\in [0,t]$ such that
\begin{equation}
|\dot{\gm}_R(s_R)|\leq D,
\end{equation}where $D$ denotes a constant independent of $R$. Let $u_R(s):=h_R(\gm_R(s),s)$. Based on the compactness of $M$, there exists a sequence $R_n$ (extracting a subsequence if necessary) such that as $n\rightarrow\infty$,
\begin{equation}
(s_{R_n},\gm_{R_n}(s_{R_n}),u_{R_n}(s_{R_n}),\dot{\gm}_{R_n}(s_{R_n}))\rightarrow (s_\infty, x_\infty,u_\infty,v_\infty).
\end{equation}
According to (L3), we have
\begin{equation}
(x_0,u_0,\dot{\gm}_\infty(0))=\Phi_{-s_\infty} (x_\infty,u_\infty,v_\infty),
\end{equation}where  $\Phi$ denotes the phase flow generated by $L$, which is conjugated to the Hamiltonian flow generated by (\ref{hjech}) via Legendre transformation. Similarly, we have
\begin{equation}
(x_0,u_0,\dot{\gm}_{R_n}(0))=\Phi_{-s_{R_n}}^{R_n} (\gm_{R_n}(s_{R_n}),u_{R_n}(s_{R_n}),\dot{\gm}_{R_n}(s_{R_n})).
\end{equation}From the construction of $\bar{L}_R$, it follows that for a given $(x,u,\dot{x})$ and any $s\in [0,t]$, as $n\rightarrow\infty$,
\begin{equation}
\Phi^{R_n}_s(x,u,\dot{x})\rightarrow \Phi_s(x,u,\dot{x}),
\end{equation}
which yields that as $n\rightarrow\infty$,
\begin{equation}
(x_0,\dot{\gm}_{{R_n}}(0))\rightarrow (x_0,\dot{x}_\infty(0)).
\end{equation}
Since $D$ is a constant independent of $R$, it follows from (L3) that there exists a constant $D'$ independent of $R$ such that
 \[\dot{\gm}_{R_n}(0)\leq D',\]
 which is in contradiction with $|\dot{\gm}_{R_n}(0)|\geq n$ for $n$ large enough.
This completes the proof of Lemma \ref{xr}.
\End

By Lemma \ref{xr}, it follows from (L3) that there exists a positive constant  $K^*$ independent of $R$ such that for any $R$ and $s\in [0,t]$, we have
\begin{equation}\label{kstar}
|\dot{\gm}_R(s)|\leq K^*.
\end{equation}
Combing Lemma \ref{A},  for a given $(x,t) (x\neq x_0, t>0)$, there exists a compact set $\Lambda_t$  independent of $R$ such that for any $R$ and $s\in [0,t]$, we have
\[(\gm_R(s),h_R(\gm_R(s),s),\dot{\gm}_R(s))\in \Lambda_t,\]
where $\gm_R:[0,t]\rightarrow M$ is a calibrated curve of $h_R$ satisfying $\gm_R(0)=x_0$ and $\gm_R(t)=x$.

\subsection{Proof of Theorem \ref{two}}

In this step, we will prove Theorem \ref{two} under the assumptions (L1)-(L4). As a preliminary, we have the following lemma:

\begin{Lemma}\label{heh} For a given $(x,t)\in M\times (0,T]$, there exists $R^*$ such that for any $R_1, R_2>R^*$,
\begin{equation}
h_{R_1}(x,t)=h_{R_2}(x,t)=h_{R^*}(x,t).
\end{equation}
\end{Lemma}
\Proof
 Let $(X_R(t),U_R(t),P_R(t))$ denote a characteristic curve generated by $H_R$, where $H_R$ denotes the Legendre transformation of $L_R$ denoted by (\ref{lr}). For the simplicity of notations, we denote
 \[\inf_{{\mathcal{S}}_{R}}U_{R}(t):=\left\{U(t):(X_R(t),U_R(t),P_R(t))\in \mathcal{S}_R\right\},\]
 where we drop the subscripts $x_0,u_0,x$ of ${\mathcal{S}}$.
  According to Theorem \ref{two1}, we have
\begin{equation}\label{huiii}
h_{R_2}(x,t)=\inf_{{\mathcal{S}}_{R_2}}U_{R_2}(t),
\end{equation}where ${\mathcal{S}}_{R_2}$ denotes the set of $(X_{R_2}(t),U_{R_2}(t),P_{R_2}(t))$ generated by $H_{R_2}$ with $X_{R_2}(0)=x_0$, $X_{R_2}(t)=x$ and $U_{R_2}(0)=u_0$.

By virtue of  (\ref{kstar}), it follows from the Legendre transformation that for any $R_1, R_2>R^*$,
\begin{equation}
(X_{R_1}(t),U_{R_1}(t),P_{R_1}(t))\in {{\mathcal{S}}_{R^*}}.
\end{equation}Hence, we have
 \[ {{\mathcal{S}}_{R_1}}\subset {{\mathcal{S}}_{R^*}}.\]
  On the other hand, it is easy to see that
\begin{equation}
 {{\mathcal{S}}_{R^*}}\subset {{\mathcal{S}}_{R_1}},
\end{equation}
Hence, we have
\begin{equation}
 {{\mathcal{S}}_{R_1}}= {{\mathcal{S}}_{R_2}}={{\mathcal{S}}_{R^*}},
\end{equation}
which together with Theorem \ref{two1} implies
\begin{equation}
h_{R_1}(x,t)=\inf_{{{\mathcal{S}}_{R_1}}}U_{R_1}(t)= \inf_{{{\mathcal{S}}_{R^*}}}U_{R^*}(t)=\inf_{{{\mathcal{S}}_{R_2}}}U_{R_2}(t)=h_{R_2}(x,t).
\end{equation}Therefore, it suffices for proving Lemma \ref{heh} to take $R^*=\max\{A^*,K^*\}$.
\End
Lemma \ref{heh} implies the existence and uniqueness of $h_{x_0,u_0}(x,t)$. Indeed, for a given $(x,t)$, one can denote
\begin{equation*}
h_{x_0,u_0}(x,t):=\lim_{R\rightarrow\infty}h_R(x,t)=h_{R^*}(x,t).
\end{equation*}In terms of Lemma \ref{heh}, we have
\begin{equation}\label{hxuxt}
h_{x_0,u_0}(x,t)=u_0+\inf_{\substack{\gm(t)=x \\  \gm(0)=y} }\int_0^tL(\gm(\tau),h_{x_0,u_0}(\gm(\tau),\tau),\dot{\gm}(\tau))d\tau,
\end{equation}where the infimums are taken among the continuous and piecewise $C^1$ curves. According to (\ref{kstar}), the uniqueness of $h_{x_0,u_0}(x,t)$ follows from a similar argument as the one in the proof of Theorem \ref{two1}. Moreover, the infimums of (\ref{hxuxt}) can be attained at a $C^1$ characteristic curve.

It is easy to see that
\begin{equation*}
\lim_{R\rightarrow\infty}{\mathcal{S}}_{R}=\mathcal{S}={\mathcal{S}}_{R^*},
\end{equation*}which implies
\begin{equation*}
\lim_{R\rightarrow\infty}\inf_{{{\mathcal{S}}_{R}}}U_{R}(t)=\inf_{{{\mathcal{S}}}}U(t).
\end{equation*}Therefore, we have

\begin{equation*}
h_{x_0,u_0}(x,t)=\inf_{{{\mathcal{S}}}}U(t)=\inf_{y\in M}\inf\left\{U(t):(X(t),U(t),P(t))\in \mathcal{S}_{y,\phi(y)}^x\right\},
\end{equation*}
where $\mathcal{S}_{y,\phi(y)}^x$ denotes the set of characteristics $(X(t),U(t),P(t))$ satisfying $X(0)=y$, $X(t)=x$ and $U(0)=\phi(y)$.

So far, we complete the proof of Theorem \ref{two} under the assumptions (L1)-(L4).

\section{\sc Representation of the viscosity solution}
In this section, we are devoted to proving Theorem \ref{one}. First of all, we construct a viscosity solution  of (\ref{hje}) with initial condition.
Based on Theorem \ref{two}, it follows that under the assumptions (L1)-(L4), there exists a unique $h_{y,\phi(y)}(x,t)\in C(M\times (0,T],\R)$ such that
\begin{equation}\label{hphi}
h_{y,\phi(y)}(x,t)=\phi(y)+\inf_{\substack{\gm(t)=x \\  \gm(0)=y} }\int_0^tL(\gm(\tau),h_{y,\phi(y)}(\gm(\tau),\tau),\dot{\gm}(\tau))d\tau,
\end{equation}where the infimums are taken among the continuous and piecewise $C^1$ curves.

A similar argument as the one in \cite{WY1} implies the following lemma.
\begin{Lemma}\label{uh}
Let
\begin{equation}
u(x,t):=\inf_{y\in M}h_{y,\phi(y)}(x,t),
\end{equation}
then \begin{equation}\label{fixu}
u(x,t)=\inf_{\gm(t)=x}\left\{\phi(\gm(0))+\int_0^tL(\gm(\tau),u(\gm(\tau),\tau),\dot{\gm}(\tau))d\tau\right\}.
\end{equation}
Moreover,
$u(x,t)$ determined by (\ref{fixu}) is a viscosity solution of (\ref{hje}) with initial condition.
\end{Lemma}

In the following, we will prove  the uniqueness of  viscosity solutions of  (\ref{hje}) with initial condition under the assumptions (H1)-(H4). More precisely, we will prove the following lemma:
\begin{Lemma}\label{main}
Under the assumptions (H1)-(H4), the viscosity solution of (\ref{hje}) with initial condition is unique.
\end{Lemma}

In order to verify Lemma \ref{main}, we need to prove the following lemmas.

\begin{Lemma}\label{lipvis}
For any $\delta>0$, a viscosity solution $u(x,t)$ of (\ref{hje}) with initial condition is necessarily Lipschitz on $M\times [\delta,T]$, and therefore satisfies (\ref{hje}) almost everywhere.
\end{Lemma}
\Proof
Let $\bar{u}(x,t)$ be a viscosity solution of (\ref{hje}) with initial condition. We consider the following equation:
\begin{equation}\label{hu}
\begin{cases}
\partial_tu(x,t)+H(x,\bar{u}(x,t),\partial_xu(x,t))=0,\\
u(x,0)=\phi(x),
\end{cases}
\end{equation}where $\bar{u}(x,t)$ is fixed. More precisely, (\ref{hu}) is equivalent to
\begin{equation}\label{hu1}
\begin{cases}
\partial_tu(x,t)+\bar{H}(x,\partial_xu(x,t))=0,\\
u(x,0)=\phi(x),
\end{cases}
\end{equation}where $\bar{H}(x,p)=H(x,\bar{u},p)$. Let $\{u_n(x,t)\}_{n\in \N}$ be a sequence of $C^2$ functions such that
\begin{equation}
\|u_n(x,t)-\bar{u}(x,t)\|_{C^0}\rightarrow 0\quad\text{as}\quad n\rightarrow \infty.
\end{equation}Let
\begin{equation}
v_n(x,t):=\inf_{\gm(t)=x}\left\{\phi(\gm(0))+\int_0^tL(\gm(\tau),u_n(\gm(\tau),\tau),\dot{\gm}(\tau))d\tau\right\}.
\end{equation} From \cite{F3}, it follows that $v_n$ is a viscosity solution of the following equation:
\begin{equation}
\begin{cases}
\partial_tu(x,t)+H_n(x,\partial_xu(x,t))=0,\\
u(x,0)=\phi(x),
\end{cases}
\end{equation}where $H_n(x,p)=H(x,u_n,p)$.
 By virtue of Proposition 4.6.6 in \cite{F3}, it follows that for each $\delta>0$, $v_n(x,t)$ is equi-Lipschitz on $M\times[\delta,T]$ where $T$ is a positive constant. Owing to the compactness of $M$,  extracting a subsequence if necessary, we have for $t\in [\delta,T]$, there exists a Lipschitz function $v(x,t)$ such that
\begin{equation}
\|v_n(x,t)-v(x,t)\|_{C^0}\rightarrow 0\quad\text{as}\quad n\rightarrow \infty.
\end{equation}.

 From the stability of viscosity solution (see Theorem 8.1.1 in \cite{F3}), it follows that $v(x,t)$ is a viscosity solution of (\ref{hu1}) on $M\times [\delta,T]$. Therefore, the comparison theorem (see \cite{Ba3} for instance) holds for the following equation:
 \begin{equation}
\partial_tu(x,t)+H(x,\bar{u}(x,t),\partial_xu(x,t))=0,
\end{equation}where $(x,t)\in M\times[\delta,T]$ and $\bar{u}$ is fixed.
Hence, we have
 \begin{equation}\label{compa}
\sup_{M\times[\delta, T]}(v(x,t)-\bar{u}(x,t))\leq\sup_M(v(x,\delta)-\bar{u}(x,\delta)).
\end{equation}
From the continuity of $v$ and $\bar{u}$, it follows that
\[\lim_{\delta\rightarrow 0}v(x,\delta)-\bar{u}(x,\delta)=v(x,0)-\bar{u}(x,0)=0\]
which together with (\ref{compa}) yields
\[v(x,t)\leq \bar{u}(x,t).\]By exchanging $v$ and $\bar{u}$ in (\ref{compa}), the comparison theorem implies
\[\bar{u}(x,t)\leq v(x,t).\]Hence,  $v(x,t)=\bar{u}(x,t)$, which means every viscosity solution of (\ref{hje}) with initial condition is Lipschitz on $M\times[\delta,T]$ for any $\delta>0$.  This completes the proof of Lemma \ref{lipvis}.
\End
Combining with Proposition \ref{semi}(b),  Lemma \ref{lipvis} implies  that every viscosity solution $u(x,t)$ of (\ref{hje}) is semiconcave on $M\times[\delta,T]$.

\textbf{Proof of Lemma \ref{main}:}
 Let \begin{equation}\label{199h}
  V(x,u,p):=H(x,u,p)-H(x,0,p).
   \end{equation}Choose a $C^2$ function $\rho(u)$ such that $\rho_R(u)=1$ for $|u|\leq R$, $\rho_R(u)=0$ for $|u|>R+1$, otherwise $0<\rho_R(u)<1$. Without loss of generality, one can require $|\rho'_R(u)|<2$. Moreover,
we denote
\begin{equation}\label{200h}
H_R(x,u,p):=H(x,0,p)+\rho_R(u)V(x,u,p).
\end{equation}
Let
\begin{equation}
\bar{H}_{R}(x,u,p):=\alpha_R (p)H_R(x,u,p)+\beta_R(p) (p^2-R^2)^2,
\end{equation}where  $\alpha_R (p)$ is a $C^2$ function satisfying
\begin{equation}
\alpha_R (p)=\left\{\begin{array}{ll}
\hspace{-0.4em}1,&  |p|\leq R+1,\\
\hspace{-0.4em}0,&|p|>R+2,\\
\end{array}\right.
\end{equation}
$\beta_R(p) $ is defined as \begin{equation}
\beta_R(p)=\left\{\begin{array}{ll}
\hspace{-0.4em}0,&  |p|\leq R,\\
\hspace{-0.4em}\mu_R,&|p|>R,\\
\end{array}\right.
\end{equation} where $\mu_R$ is a sufficient large constant.

 $\bar{H}_{R}(x,u,p)$  generates the following equation:
 \begin{equation}\label{Hr}
\partial_tu(x,t)+\bar{H}_{R}(x,u(x,t),\partial_xu(x,t))=0.
\end{equation}
It is easy to see that $\bar{H}_{R}(x,u,p)$ is uniformly Lipschitz with respect to $u$, for which the comparison theorem holds (see \cite{Ba3}).

Let $u(x,t)$ be a viscosity solution of (\ref{hje}). Let $\mathcal{D}$ be  the differentiable points of $u(x,t)$ on $M\times[\delta,T]$. It follows from Lemma \ref{lipvis} that $\mathcal{D}$ has full Lebesgue measure. For $(x,t)\in \mathcal{D}$, we have
\begin{equation}
|\partial_xu(x,t)|\leq K_1,
\end{equation}where $K_1$ is a constant independent of $(x,t)$.
Based on the compactness of $M\times[\delta,T]$, it yields that $u(x,t)$ has a bound denoted by $K_2$. Taking $R=\max\{K_1,K_2\}$, it follows from the construction of $\bar{H}_{R}(x,u,p)$ that for $(x,t)\in \mathcal{D}$, we have
\begin{equation}\label{urb}
\partial_tu(x,t)+H_{\tilde{R}}(x,u(x,t),\partial_xu(x,t))=0.
\end{equation}
By Proposition \ref{semi}(a), $u(x,t)$ is a viscosity solution of (\ref{Hr}) on $M\times[\delta,T]$. Hence, the comparison theorem still holds for the viscosity solutions of (\ref{hje}) on $M\times[\delta,T]$.
 Let $u_1(x,t)$ and $u_2(x,t)$ be two viscosity solutions of (\ref{hje}) with initial condition. By the comparison theorem, it follows that
 \begin{equation}
\sup_{M\times[\delta, T]}(u_1(x,t)-u_2(x,t))\leq\sup_M(u_1(x,\delta)-u_2(x,\delta)),
\end{equation}which together  the continuity of $u_1$ and $u_2$ yields
 \[u_1(x,t)=u_2(x,t).\]
Therefore, we obtain the uniqueness of the viscosity solution of  (\ref{hje}) with initial condition.\End

So far,  we complete the proof of Theorem \ref{one}.

 \vspace{2ex}
\noindent\textbf{Acknowledgement}
This work is partially under the support of National Natural Science Foundation of China (Grant No. 11171071,  11325103) and
National Basic Research Program of China (Grant No. 11171146).

\vspace{2em}

{\sc Lin Wang}

{\sc School of Mathematical Sciences, Fudan University,
Shanghai 200433,
China.}

 {\it E-mail address:} \texttt{linwang.math@gmail.com}

\vspace{1em}

{\sc Jun Yan}

{\sc School of Mathematical Sciences, Fudan University,
Shanghai 200433,
China.}

 {\it E-mail address:} \texttt{yanjun@fudan.edu.cn}


\addcontentsline{toc}{section}{\sc References}

\begin{thebibliography}{999}
\renewcommand{\baselinestretch}{1.0}
\setlength\itemsep{-1pt}
\small


\bibitem{Ar} V. I. Arnold.
{\it Geometric methods in the theory of ordinary differential
equations}. Springer-Verlag, New York. 1983.


\bibitem{Au} S. Aubry. {\it The twist map, the extended Frenkel-Kontorova model and the devil's staircase}. Phys. D. \textbf{7} (1983), 240-258.


\bibitem{AD} S. Aubry and P. Y. Le Daeron. {\it The discrete Frenkel-Kontorova model and its extensions I: exact results for the ground
 states}. Phys. Rev. D \textbf{8} (1983), 381-422.



\bibitem{Ba1} G. Barles. {\it Existence results for first order Hamilton-Jacobi equations}. Ann. Inst. Henri. Poincar\'{e}, \textbf{1}(5):325-340, 1984.

\bibitem{Ba2} G. Barles. {\it Remarques sur des r\'{e}sultants d'existence pour les \'{e}quations de Hamilton-Jacobi du premier ordre}. Ann. Inst. Henri. Poincar\'{e}, \textbf{2}(1):21-32, 1985.



\bibitem{BaS} G. Barles and P. E. Souganidis. {\it On the large time behavior of solutions of Hamilton-Jacobi equations}. SIAM J. Math. Anal, \textbf{31}:925-939, 2000.

\bibitem{Ba3} G. Barles. {\it An introduction to the theory of viscosity solutions for first-order Hamilton-Jacobi equations and applications}. Hamilton-Jacobi Equations: Approximations, Numerical Analysis and Applications, Lecture Notes in Mathematics \textbf{2074} Springer-Verlag Berlin Heidelberg 2013.



\bibitem{Be} P. Bernard. {\it Existence of $C^{1,1}$ critical sub-solutions of the Hamilton-Jacobi equation on compact manifolds}. Annales Scientifiques de l'\'{E}cole Normale Sup\'{e}rieure, \textbf{40}(3):445-452, 2007.




\bibitem{CS}   P. Cannarsa and C. Sinestrari. {\it Semiconcave functions, Hamilton-Jacobi equations, and optimal control}. Vol. \textbf{58}. Springer, 2004.

\bibitem{CIPP} G. Contreras, R. Iturriaga, G. P. Paternain and M. Paternain.
{\it Lagrangian graphs, minimizing measures and Ma\~{n}\'{e}'s
critical values}. Geom. Funct. Anal. \textbf{8} (1998), 788-809.



\bibitem{CEL} M. G. Crandall, L. C. Evans, and P.-L. Lions. {\it Some properties of viscosity solutions
of Hamilton-Jacobi equations}. Trans. Amer. Math. Soc., \textbf{282}(2):487-502,
1984.

\bibitem{CHL1} M. G. Crandall, H. Ishii, and P.-L. Lions. {\it Uniqueness of viscosity solutions of Hamilton-Jacobi equations revisited}. J. Math. Soc. Japan, \textbf{39}(4):581-596, 1987.

\bibitem{CHL2} M. G. Crandall, H. Ishii, and P.-L. Lions. {\it User¡¯s guide to viscosity
solutions of second order partial differential equations}. Bull. Amer. Math.
Soc. (N.S.), \textbf{27}(1):1-67, 1992.



\bibitem{CL1}M. G. Crandall and P.-L. Lions. {\it Condition d¡¯unicit\'{e} pour les solutions
g\'{e}n\'{e}ralis\'{e}es des \'{e}quations de Hamilton-Jacobi du premier ordre}. C. R.
Acad. Sci. Paris S\'{e}r. I Math., \textbf{292}(3):183-186, 1981.

\bibitem{CL2}M. G. Crandall and P.-L. Lions. {\it Viscosity solutions of Hamilton-
Jacobi equations}. Trans. Amer. Math. Soc., \textbf{277}(1):1-42, 1983.
%





\bibitem{E}W. E. {\it Aubry-Mather theory and periodic solutions of the forced Burgers
equation}. Comm. Pure Appl. Math., \textbf{52}(7):811-828, 1999.


\bibitem{F1} A. Fathi. {\it Solutions KAM faibles conjugu\'{e}es et barri\`{e}res de Peierls}. C. R.
Acad. Sci. Paris S\'{e}r. I Math., \textbf{325}(6):649-652, 1997.

 \bibitem{F2} A. Fathi. {\it Th\'{e}or\`{e}me KAM faible et th\'{e}orie de Mather sur les syst\`{e}mes lagrangiens}.
C. R. Acad. Sci. Paris S\'{e}r. I Math., \textbf{324}(9):1043-1046, 1997.


 \bibitem{F22} A. Fathi. {\it Sur la convergence du semi-groupe de Lax-Oleinik}.
C. R. Acad. Sci. Paris S\'{e}r. I Math., \textbf{327}(3):267-270, 1998.




 \bibitem{F3} A. Fathi. {\it Weak KAM Theorem in Lagrangian Dynamics}. Preliminary Version Number 10, 2008.

 \bibitem{FS} A. Fathi and A. Siconolfi. {\it Existence of $C^1$ critical subsolutions of the
Hamilton-Jacobi equation}.
Invent. math., \textbf{155}:363-388, 2004.



\bibitem{GT}   D. Gilbarg and N. Trudinger. {\it Elliptic partial differential equations of second order}. Originally published as volume \textbf{224} in the series: Grundlehren der mathematischen Wissenschaften. Reprint of the 1998 ed. Springer-Verlag 2003 .


 \bibitem{HS} P.-F. Hsieh and Y. Sibuya. {\it Basic theory of ordinary differential equations}.
Springer-Verlag New York 1999.


\bibitem{I1}H. Ishii. {\it Existence and uniqueness of solutions of Hamilton-Jacobi equations}. Funkcialaj Ekvacioj, \textbf{29}:167-188, 1986.




\bibitem{Li} P.-L. Lions. {\it Generalized solutions of Hamilton-Jacobi equations}. volume
69 of Research Notes in Mathematics. Pitman (Advanced Publishing Program),
Boston, Mass., 1982.

 \bibitem{M0} J. N. Mather. {\it Existence of quasi periodic orbits for twist homeomorphisms of the
 annulus}. Topology \textbf{21} (1982), 457-467.



\bibitem{Mat} J. N. Mather. {\it More Denjoy minimal sets for area preserving diffeomorphisms}. Comment. Math. Helv. \textbf{60} (1985), 508-557.

 \bibitem{M2} J. N. Mather. {\it A criterion for the non-existence of
 invariant circle}. Publ. Math. IHES \textbf{63} (1986), 301-309.

 \bibitem{M3} J. N. Mather. {\it Modulus of continuity for Peierls's
 barrier}. Periodic Solutions of Hamiltonian Systems and Related
 Topics. ed. P.H.Rabinowitz {\it et al}. NATO ASI Series C \textbf{209}.
 Reidel: Dordrecht, (1987), 177-202.



\bibitem{M1}J. N. Mather. {\it Action minimizing invariant measures for positive definite Lagrangian
systems}. Math. Z., \textbf{207}(2):169-207, 1991.

%
\bibitem{Mvc} J. N. Mather. {\it Variational construction of connecting orbits}. Ann. Inst. Fourier
(Grenoble), \textbf{43}(5):1349-1386, 1993.

%



\bibitem{os} W. F. Osgood. {\it Bewise der Existenz einer L\"{o}sung einer Differentialgleichung $dy/dx=f(x,y)$ ohne Heinzunahme der Cauchy-Lipschizschen Bedingung}. Monat. Math.-Phys., \textbf{9}:331-345, 1898.

\bibitem{SWY} X. Su, L. Wang and J. Yan. {\it Weak KAM theory for Hamilton-Jacobi equations I: the solution semigroup under proper conditions}. preprint.



\bibitem{WY1}  L. Wang and J. Yan. {\it Weak KAM theory for Hamilton-Jacobi equations II: the fundamental solution under  Lipschitz conditions}. preprint.

\end{thebibliography}
\end{document}